\documentclass{amsart}
\newcommand{\note}{\noindent {\bf Notation. }}
\newcommand{\prop}{\noindent {\bf Properties. }}

\newcommand{\ws}{\hspace{4pt}}
\newtheorem{theorem}{Theorem}
\newtheorem{statement}{Statement}

\newtheorem{proposition}{Proposition}
\newtheorem{cor}{Corollary}
\newtheorem{lemma}{Lemma}

\newtheorem{defi}{Definition}
\usepackage[]{amssymb,amsopn}
\usepackage[]{amsmath}
\usepackage[]{eucal}
\usepackage[]{latexsym}
\usepackage{color}

\makeatletter
\@namedef{subjclassname@2020}{%
  \textup{2020} Mathematics Subject Classification}
\makeatother

\begin{document}

\title[Fractional Bessel-Sobolev spaces]{Fractional Bessel-Sobolev and Bessel B-potential spaces}
\author{Mouna Chegaar and \'A. P. Horv\'ath}

\thanks{}

\begin{abstract} After proving the equivalence of the Bessel $K$-functional and the corresponding spherical modulus of smoothness we define fractional Bessel-Sobolev spaces. As an analog of the classical one the imbedding relation of fractional Bessel-Sobolev and Bessel B-potential spaces is pointed out. Applying the defined notions a potential theoretic characterization of removable sets with respect to certain fractional partial differential operator is derived.

\end{abstract}
\maketitle

\section{Introduction}

Bessel and fractional Bessel differential operators occur frequently in physics and in practice, see e.g. \cite{mg} or \cite{esk}. The generalized (or Bessel) translation defined by Delsarte, \cite{d} and developed by Levitan, \cite{le} has proven to be a fairly effective tool for examining these operators. Application  of Bessel translation rather than the standard one modifies the geometry of the space appropriately. This is why it is worth to introduce B-Hardy, B-Morrey or the different generalized Bessel B-potential spaces  etc., see \cite{hs}, \cite{gh}, \cite{sk} and \cite{esk2}.  \\
For instance, to study the convergence properties of a numerical solution, to characterize removable sets or to develop the suitable spectral theory the investigation of the underlying space is essential.\\
First we focus on the equivalence of modulus of smoothness and the corresponding $K$-functional. It is a crucial inequality in approximation theory and also in the theory of interpolation of spaces. In \cite{gott} several results on $K$-functional are summarized in classical cases. In the following we extend the one-dimensional results of \cite{p} to higher dimensions. The main tool of this section is the weighted spherical modulus of smoothness introduced in \cite{ch}.\\
This equivalence allows to define the appropriate spherical Gagliardo-type seminorm (c.f. \cite{dpv}) and fractional Bessel-Sobolev spaces, i.e. to extend the results of \cite{p} and \cite{sk} to fractional exponents.\\
As usual, Bessel translation defines a convolution and so Bessel B-potential spaces, i.e. the spaces of functions derived by Bessel convolution with the Bessel kernel. Classical Bessel potential spaces were introduced systematically in a series of papers by Aronszajn and Smith, see e.g. \cite{as} and by Calder\'on, see \cite{c} for the one-dimensional case. An analog of Calder\'on's theorem is derived in the Bessel case, see \cite{sk} or \cite{h}. The fractional case is a little different, c.f. \cite{a}. By examining the defined norms, a suitable imbedding theorem can be proven.\\
As an application, we characterize removable sets with respect to certain fractional partial differential operator. The potential theoretic characterization of removable sets is due to Carleson, see \cite{ca}. For Laplace-Bessel type operators similar results can be found in \cite{h}. For fractional parabolic operators the results of \cite{mp} are extended in \cite{ch}. In the last section we give a characterization of removable sets with respect to fractional elliptic operators via capacity defined by fractional Bessel B-potencial spaces.

\section{Notation, preliminaries}

\subsection{Norms, spaces}
On $\mathbb{R}^{n}_{+}:=\{(x_1, \dots ,x_n): x_i>0 \ws i=1, \dots ,n\}$ we introduce the $L^p_a(\mathbb{R}^{n}_{+})$ space with norm
$$\|f\|_{p,a}^p=\int_{\mathbb{R}^n_{+}}|f(x)|^px^adx,$$
where $a=(a_1, \dots , a_n)$, $a_i>0$, $i=1, \dots , n$; $x^a=\prod_{i=1}^nx_i^{a_i}$.\\
The appropriate weighted inner product is as follows.
$$\langle f,g\rangle_a:=\int_{\mathbb{R}^n_{+}}f(x)g(x)x^adx.$$
$\mathcal{S}_{(e)}$ stands for the even functions from the Schwartz class, i.e.\\ $\mathcal{S}_{(e)}:=\left\{f\in \mathcal{S} : \frac{\partial^{2m+1}}{\partial x_i^{2m+1}}f\large{|}_{x_i=0}=0, \ws \forall m \in \mathbb{N}, \ws i=1, \dots , n\right\}.$ \\
If $h\in \mathcal{S}_{(e)}'$, $\varphi \in \mathcal{S}_{(e)}$ we denote the effect of $h$ on $\varphi$ by $\langle h,\varphi\rangle_a$ as well.

\subsection{Bessel differential operators}
The Bessel-Laplace operator on $\mathbb{R}^{n}_{+}$ is defined as follows.
\begin{equation}\label{bdo}\Delta_a:=\sum_{i=1}^n B_{\alpha_i}, \ws \mbox{where} \ws B_{\alpha_i}:=\frac{\partial^2}{ \partial x_i\partial x_i}+\frac{2\alpha_i+1}{x_i}\frac{\partial}{\partial x_i},\end{equation}
and
$$a_i=2\alpha_i+1, $$
(i.e. $\alpha_i>-\frac{1}{2}$, $i=1, \dots , n$).
Hereafter we deal with partial differential operators which contain integer or fractional powers of $\Delta_a$ or $I-\Delta_a$.

\subsection{Bessel functions, Bessel translation} The entire Bessel functions are
\begin{equation}\label{B}j_\alpha(z)=\Gamma(\alpha+1)\left(\frac{2}{z}\right)^\alpha J_\alpha(z)=\sum_{k=0}^\infty\frac{(-1)^k\Gamma(\alpha+1)}{\Gamma(k+1)\Gamma(k+\alpha+1)}\left(\frac{z}{2}\right)^{2k},\end{equation}
where $\alpha>-\frac{1}{2}$. We note that $j_{-\frac{1}{2}}(x)=\cos x$. $j_\alpha(\lambda z)$ are the eigenfunctions of $B_\alpha$.
\begin{equation}\label{j1}\|j_\alpha\|_{\infty, \mathbb{R}_+}=j_\alpha(0)=1.\end{equation}
\begin{equation}\label{j2}j_\alpha'(z)=-\frac{1}{2(\alpha+1)}zj_{\alpha+1}(z),\end{equation}
see \cite{be}.

We introduce the following abbreviation, cf. \cite{ss}.
$$\mathbb{j}_a(x,\xi):=\prod_{i=1}^nj_{\alpha_i}(x_i\xi_i).$$

Subsequently we also need the modified Bessel function of the second kind $K_\nu$.
\begin{equation}\label{k0} K_\nu(x)=\frac{\pi}{2}\frac{i^{\nu}J_{-\nu}(ix)-i^{-\nu}J_\nu(ix)}{\sin\nu \pi},\end{equation}
see \cite{be}. Growth properties of $K_\nu$ are as follows.

\noindent Around zero:
\begin{equation}\label{k1}K_\nu(r)\sim \left\{\begin{array}{ll}-\ln\frac{r}{2}-c, \ws \mbox{if}\ws  \nu=0\\ C(\nu)r^{-\nu},  \ws \mbox{if}\ws \nu>0,\end{array}\right.\end{equation}

\noindent Around infinity:
\begin{equation}\label{k2}K_\nu(r)\sim \frac{c}{\sqrt{r}}e^{-r},\end{equation}
see e.g. \cite[(1.2.20)-(1.2.22)]{ah}.

\noindent Moreover around zero for $\nu\ge 0$ we have
\begin{equation}\label{k3}|(r^\nu K_\nu(r))'|\le c \left\{\begin{array}{ll}r^{2\nu-1}\ln\frac{1}{r}, \ws \mbox{if}\ws  1>2\nu, \ws \nu \in\mathbb{N}\\ r^{2\nu-1},  \ws \mbox{if}\ws  1>2\nu, \ws \nu \not\in\mathbb{N}\\
1,  \ws \mbox{if}\ws 1<2\nu, \ws \nu \not\in\mathbb{N}\\
r,  \ws \mbox{if}\ws 1<2\nu, \ws \nu \not\in\mathbb{N}
,\end{array}\right.,\end{equation}
see \cite[(40),(41)]{dls}.\\

With the above notation the Bessel translation of a function, $f$ (see e.g. \cite{le}, \cite{p}, \cite{ss}) is
$$T_{a}^tf(x)=T_{a_n}^{t_n}\dots T_{a_1}^{t_1}f(x_1,\dots ,x_n),$$
where
\begin{equation}\label{tra1}T_{a_i}^{t_i}f(x_1,\dots ,x_n)$$ $$=\frac{\Gamma(\alpha_i+1)}{\sqrt{\pi}\Gamma\left(\alpha_i+\frac{1}{2}\right)}\int_0^\pi f(x_1,\dots,\sqrt{x_i^2+t_i^2 -2x_it_i\cos\vartheta_i}, x_{i+1}, \dots, x_n)\sin\vartheta^{2\alpha_i}d\vartheta_i.\end{equation}
$$T_{a}^tf(x)=T_{a}^xf(t),$$
\begin{equation}\label{tt}T_{a}^yf(x)T_{a}^zf(x)=T_{a}^zf(x)T_{a}^yf(x),\end{equation}
see \cite[(7.1)]{le}.
For all $a$ ($|a|\ge 0$)$, T_a$ is positive operator, and
\begin{equation}\label{t1}\|T_{a}^tf(x)\|_{p,a}\le \|f\|_{p,a}, \ws \ws 1\le p \le\infty,\end{equation}
$1\le p \le\infty$, see e.g. \cite{le}.

\subsection{Bessel convolution, Hankel transformation}
The generalized convolution with respect to the Bessel translation is
$$f*_a g =\int_{\mathbb{R}^n_+}T_{a,x}^tf(x)g(x)x^adx.$$
We have
\begin{equation}\label{cca}f*_a g =g*_a f,\ws \ws f*_a(g*_ah)= (f*_ag)*_ah.\end{equation}
Furthermore Young's inequality fulfils i. e. if $1\le p, q, r \le \infty$ with $\frac{1}{r}=\frac{1}{p}+\frac{1}{q}-1$; if $f\in L^p_\alpha$ and $g \in L^q_\alpha$, then
\begin{equation}\label{Y}\|f*g\|_{r,a} \le \|f\|_{p,a}\|g\|_{q,a},\end{equation}
see \cite[(3.178)]{ss}.

$$\mathcal{H}_a(f,\xi)=\hat{f}(\xi)=\int_{\mathbb{R}^n_+}f(x)\mathbb{j}_a(x,\xi)x^adx.$$
If $f\in L^1_a(\mathbb{R}^n_+)$ $(|a|>0)$ is of bounded variation in a neighborhood of a point $x$ of continuity of $f$, then the next inversion formula holds, see [page 38.]\cite{ss}.
\begin{equation}\label{bi}f(x)=\mathcal{H}_a^{-1}(\hat{f}(\xi))(x)=\frac{2^{n-|a|}}{\prod_{i=1}^n\Gamma^2(\alpha_i+1)}\int_{\mathbb{R}^n_+}\hat{f}(\xi)\mathbb{j}_a(x,\xi)\xi^ad\xi.\end{equation}

The relationship of Hankel transform to translation and convolution is as follows.
\begin{equation}\label{HT} \mathcal{H}_a(T^y_af)(\xi)= \mathbb{j}_a(u,y)(\xi,y)\mathcal{H}_a(f)(\xi); \ws \ws \ws T^y_a(\mathcal{H}_a f)(\xi)=\mathcal{H}_a(\mathbb{j}_a(u,y)f(x)),\end{equation}
see  \cite[(3.158)]{ss}.\\
for $f,g \in \mathcal{S}_e$ we have
\begin{equation}\label{CH} \mathcal{H}_a(f*_ag)=\hat{f}\hat{g},\end{equation}
see e.g. \cite[(3.176)]{ss}.

 Subsequently we need that the Hankel transform of a radial function. We begin with the following formula.
 \begin{equation}\label{besatl} \frac{1}{|S_{1+}^n|_a}\int_{S_{1+}^n}\mathbb{j}_a(r\Theta,\xi))\Theta^adS(\Theta)=j_{\frac{n+|a|}{2}-1}(r|\xi|),\end{equation}
  c.f. \cite[(3.140)]{ss}.\\

 Then let $f$ be a radial function, i.e. $f(x)=\varphi(|x|)$, where $\varphi$ is a function of one variable. Denoting by $r=|x|$ $\mathcal{H}_a(f)$, is also radial, furthermore we have
\begin{equation}\label{rh}\mathcal{H}_a(f)(\xi)=|S_{1+}^n|_a\int_0^\infty \varphi(r)j_{\frac{n+|a|}{2}-1}(|\xi|r)r^{n+|a|-1}dr,\end{equation}
where
 \begin{equation}\label{FS}|S_{1+}^n|_a=\int_{S_{1+}}\Theta^adS(\Theta)=\frac{\prod_{i=1}^n \Gamma(\alpha_i+1)}{2^{n-1}\Gamma\left(\frac{n+|a|}{2}\right)},\end{equation}
 and $dS$ stands for the integration on the sphere; see \cite[Lemma 3.2]{esk}.

\subsection{Differences, moduli of smoothness}
\begin{defi}
 The differences are
 $$\Delta_{y,a}f(x):=f(x)-T^y_af(x); \ws \ws \Delta_{y,a}^kf(x):=\Delta_{y,a}\Delta_{y,a}^{k-1}f(x),$$
 $$\tilde{\Delta}_{y,a}^kf(x):=\sum_{l=0}^k(-1)^l \binom{k}{l} T^{\sqrt{l}y}_af(x).$$
 \end{defi}

 \note
 $$\mathcal{J}_{\alpha,k}(x):=\sum_{l=0}^k(-1)^l \binom{k}{l}j_\alpha(\sqrt{l}x).$$

 Subsequently the following property of $\mathcal{J}_{\alpha,k}$ will be useful.

 \medskip

 \begin{lemma}\cite[Lemma 4.1]{p}
 \begin{equation}\label{fo}\mathcal{H}_{\alpha}^{-1}\left(\frac{\mathcal{J}_{\alpha,k}(\xi)}{\xi^{2k}}\right)\in L^1_\alpha(\mathbb{R}_+).\end{equation}
 \end{lemma}

 \medskip

 With this notation we introduce  the following spherical differences.

  \medskip

 \begin{defi}\label{ms}
 $$\Delta_{sph,r,a}f(x):=\frac{1}{|S^n_{1+}|_a}\int_{S^n_{1+}}(f(x)-T^{r\Theta}_af(x))\Theta^adS(\Theta);$$
 $$\Delta_{sph,r,a}^kf(x):=\Delta_{sph,r,a}\Delta_{sph,r,a}^{k-1}f(x),$$
 and
  $$\tilde{\Delta}_{sph,r,a}^kf(x):=\frac{1}{|S^n_{1+}|_a}\int_{S^n_{1+}}\sum_{l=0}^k(-1)^l  \binom{k}{l} T^{\sqrt{l}r\Theta}_af(x)\Theta^adS(\Theta),$$
 c.f. \eqref{FS}.
 $$\omega_{sph,m}(f,t)_{a,p}:=\sup_{0<h<t}\|\Delta_{sph,h,a}^mf\|_{p,a}.$$
 \end{defi}

 \medskip

 \prop

 In view of \eqref{HT} and \eqref{besatl} we have
 \begin{equation}\label{spht} \mathcal {H}_a(\Delta_{sph,r,a}f)(\xi)=(1-j_{\frac{n+|a|}{2}-1}(r|\xi|))(\mathcal{H}_a f)(\xi).\end{equation}

 From the definition and from \eqref{t1} it follows that $\omega_{sph,m}(f,t)_{a,p}$ is increasing in $t$, moreover
 \begin{equation}\label{o1} \omega_{sph,m}(f+g,t)_{a,p}\le \omega_{sph,m}(f,t)_{a,p}+\omega_{sph,m}(g,t)_{a,p} \end{equation}
 and
 \begin{equation}\label{o2} \omega_{sph,m}(f,\delta)_{a,p}\le 2^m \|f\|_{p,a}. \end{equation}

 \section{K-functional}

 We give the definition of the Bessel-Sobolev space and the corresponding $K$-functional, see \cite{sk}, \cite{h} and in 1-dimension see \cite[(1.9)]{p}. The following theorem describes the connection beetween the defined $K$-functional and spherical modulus of smoothness.

 \medskip

 \begin{defi}\label{Sob} Let  $1\le p<\infty$; $m\in\mathbb{N}$.
 $$W^{m,p}_{\Delta_a}:=\left\{f: \mathbb{R}^n_{+} \to \mathbb{R} : \Delta_a^k f \in  L^p_a, \ws k=0, \dots , m \right\},$$ $$ \|f\|_{W^{m,p}_{\Delta_a}}=\left(\sum_{k=0}^m \|\Delta_a^k f\|_{p,a}^p\right)^{\frac{1}{p}}.$$
 \end{defi}

 \medskip

\begin{defi}\label{k}
$$K\left(f, t, L_a^p, W^{k,p}_{\Delta_a}\right)=\inf\{\|f-g\|_{p,a}+t\|\Delta_a^kg\|_{p,a} : g\in W^{k,p}_{\Delta_a}\},$$
 where $f\in L_a^p$, $t>0$.
\end{defi}

\medskip

\begin{theorem}\label{ok} There is a positive constant $c=c(m,n,a)$, such that
$$\frac{1}{c}\omega_{sph,m}(f,t)_{a,p}\le K\left(f, t^{2m}, L_a^p, W^{m,p}_{\Delta_a}\right)\le c \omega_{sph,m}(f,t)_{a,p},$$
 where $f\in L_a^p$, $t>0$.
\end{theorem}

\medskip

This theorem is an $n$-dimensional version of \cite[Theorem 1.2]{p}. The proof follows the chain of ideas of the original one. The main tool for the $n$-dimensional case is the spherical modulus of smoothness defined above. For sake of self-containedness we give the detailed proof.

For the proof of the inequality on the left-hand side we need the next lemma.

\medskip

\begin{lemma}\label{lbal} If $ g\in W^{k,p}_{\Delta_a}$, then for all $t>0$ and $1\le p \le \infty$
$$\omega_{sph,m}(g,t)_{a,p}\le c t^{2m} \|\Delta_a^m g\|_{p,a},$$
where $c=c(a,n,m)$.
\end{lemma}

\medskip

\proof (of the left inequality of Theorem \ref{ok}) Let $f\in L_a^p$ and $g\in W^{m,p}_{\Delta_a}$. According to \eqref{o1} and \eqref{o2} we have
$$\omega_{sph,m}(f,t)_{a,p}\le \omega_{sph,m}(f-g,t)_{a,p}+\omega_{sph,m}(g,t)_{a,p}\le 2^m\|f-g\|_{p,a}+ c t^{2m} \|\Delta_a^m g\|_{p,a}$$ $$\le c\left(\|f-g\|_{p,a}+ t^{2m} \|\Delta_a^m g\|_{p,a}\right),$$
where Lemma \ref{lbal} is applied. Taking infimum over $g\in W^{k,p}_{\Delta_a}$, the left-hand side of Theorem \ref{ok} is proved.

\medskip

To prove Lemma \ref{lbal} we need some definitions and further lemmas.

\note

First we define the radial function $\eta \in \mathcal{S}_e$ such that $|\eta|\le 1$,\\ $\eta(x)=\eta(|x|)=\left\{\begin{array}{ll}1, \ws |x|\le 1,\\0, \ws |x|\ge 2\end{array}\right.$.\\
For an $f\in L_a^p$, let
\begin{equation}\label{Pnu}P_\nu(f)(x):=\mathcal{H}_a^{-1}\left(\eta\left(\frac{|\xi|}{\nu}\right)\hat{f}(\xi)\right)=\varrho_\nu*_af(x),\end{equation}
where $\varrho_\nu(x):=\mathcal{H}_a^{-1}\left(\eta\left(\frac{|\xi|}{\nu}\right)\right)$.\\
Let us denote by $E_e$ the set of even entire functions. We also define the space
$$M(\nu, p,a):=\{h\in L^p_a\cap E_e: \mathrm{supp}\hat{h}\subset B(0,\nu)\}.$$

\prop

\begin{equation}\label{pnorm}\|P_\nu(f)\|_{p,a}\le c \|f\|_{p,a}; \ws \ws c\neq c(\nu).\end{equation}
Indeed, in view of \eqref{Y} $\|P_\nu(f)\|_{p,a}\le \|\varrho_\nu\|_{1,a}\|f\|_{p,a}$. Since  $\eta \in \mathcal{S}_e$,  $\hat{\eta} \in \mathcal{S}_e\subset L^1_a$. Finally it can be readily seen, that $\|\varrho_\nu\|_{1,a}=\|\varrho_1\|_{1,a}$.

\begin{equation}\label{Psupp}\mathrm{supp}\hat{P_\nu}(f)\subset B(0,\nu).\end{equation}
Indeed, in distribution sense $\hat{P_\nu}(f)=\eta\left(\frac{|\xi|}{\nu}\right)\hat{f}(\xi)$, cf. \cite[Lemma 4.5]{p}.

Since $\eta \in \mathcal{S}_e$, \eqref{pnorm} and \eqref{Psupp} ensure that $P_\nu(f)\in M(\nu, p,a)$.

\begin{equation}\label{Pg}P_\nu(h)=h, \ws \ws \forall \ws h \in M(\nu, p,a). \end{equation}
Indeed, if $\mathrm{supp}\hat{h}\subset B(0,\nu)$, $\hat{h}(\xi)\eta\left(\frac{|\xi|}{\nu}\right)=\hat{h}(\xi)$.

Moreover we define the best approximation of an $f\in L_a^p$ with functions from $M(\nu, p,a)$.

\medskip

\begin{defi}
$$E_\nu(f)_{p,a}:=\inf\{\|f-h\|_{p,a} : h \in M(\nu, p,a)\}.$$
\end{defi}

\medskip

A further important property of $P_\nu(f)$ is as follows.

\begin{equation}\label{E} \|f-P_\nu(f)\|_{p,a}\le c E_\nu(f)_{p,a}; \ws \ws f\in L_a^p, \ws c\neq c(f).\end{equation}
Indeed, let $h\in  M(\nu, p,a)$ such that $\|f-h\|_{p,a}\le 2 E_\nu(f)_{p,a}$. Then by \eqref{Pg} and \eqref{pnorm} we have $\|f-P_\nu(f)\|_{p,a}\le \|f-h\|_{p,a}+\|P_\nu(f-h)\|_{p,a}\le cE_\nu(f)_{p,a}$.

\medskip

\begin{lemma}\label{l2}Let  $1\le p \le \infty$ and $g\in W^{m,p}_{\Delta_a}$. Then we have
$$\|\Delta_{sph,h,a}^mP_\nu(g)\|_{p,a}\le c h^{2m}\|\Delta_a^m g\|_{p,a},$$
where $c=c(n,a,m)$.
\end{lemma}

\proof According to \eqref{spht} we have
$$\mathcal{H}_a\left(\Delta_{sph,h,a}^mP_\nu(g)\right)(\xi)=\left(1-j_{\frac{n+|a|}{2}-1}(h|\xi|\right)^m\eta\left(\frac{|\xi|}{\nu}\right)\hat{g}(\xi)$$
$$= h^{2m}\eta\left(\frac{|\xi|}{\nu}\right)\frac{\left(1-j_{\frac{n+|a|}{2}-1}(h|\xi|)\right)^m}{(-1)^m(h|\xi|)^{2m}}(-1)^m|\xi|^{2m}\hat{g}(\xi).$$
According to \eqref{fo} $u_h:=\mathcal{H}_a^{-1}\left(\frac{\left(1-j_{\frac{n+|a|}{2}-1}(h|\xi|)\right)}{(h|\xi|)^2}\right)\in L^1_a$ and its norm is independent of $h$, so by \eqref{Y} the same is valid for $u_h*_a\dots *_a u_h=:\psi_h$. That is by \eqref{Y} and \eqref{CH} we have
$$\|\Delta_{sph,h,a}^mP_\nu(g)\|_{p,a}\le h^{2m}\|\varrho_\nu\|_{1,a}\|\psi_h\|_{1,a}\|\Delta_a^m g\|_{p,a}\le  c h^{2m}\|\Delta_a^m g\|_{p,a},$$
where $c=\|\varrho_1\|_{1,a}\|\psi_1\|_{1,a}$.

\medskip

\begin{lemma}\label{l3}Let $1\le p \le \infty$ and $g\in W^{m,p}_{\Delta_a}$. Then we have
$$\|\tilde{\Delta}_{sph,h,a}^mg\|_{p,a}\le c h^{2m}\|\Delta_a^m g\|_{p,a},$$
where $c=c(n,a,m)$.
\end{lemma}

\proof
By Fubini's theorem, \eqref{HT} and \eqref{besatl} we have
$$\mathcal{H}_a\left(\tilde{\Delta}_{sph,h,a}^mg\right)(\xi)=\hat{g}(\xi)\mathcal{J}_{\frac{n+|a|}{2}-1,m}(h|\xi|).$$
In view of \eqref{fo} $\mathcal{H}_a^{-1}\left(\frac{\mathcal{J}_{\frac{n+|a|}{2}-1,m}(h|\xi|)}{|\xi|^{2m}}\right)\in L^1_a$, and $\left\|\mathcal{H}_a^{-1}\left(\frac{\mathcal{J}_{\frac{n+|a|}{2}-1,m}(h|\xi|)}{|\xi|^{2m}}\right)\right\|_{1,a}=h^{2m}\left\|\mathcal{H}_a^{-1}\left(\frac{\mathcal{J}_{\frac{n+|a|}{2}-1,m}(|\xi|)}{|\xi|^{2m}}\right)\right\|_{1,a}$. Thus by \eqref{Y} we have
$$\|\tilde{\Delta}_{sph,h,a}^mg\|_{p,a}=\left\|\mathcal{H}_a^{-1}\left(\frac{\mathcal{J}_{\frac{n+|a|}{2}-1,m}(h|\xi|)}{|\xi|^{2m}}\right)*_a\Delta_a^m g\right\|_{p,a}\le c h^{2m}\|\Delta_a^m g\|_{p,a},$$
where $c=\left\|\mathcal{H}_a^{-1}\left(\frac{\mathcal{J}_{\frac{n+|a|}{2}-1,m}(|\xi|)}{|\xi|^{2m}}\right)\right\|_{1,a}$.

\medskip

\begin{lemma}\label{l4} Let  $1\le p \le \infty$ and $g\in W^{m,p}_{\Delta_a}$. Then we have
$$E_\nu(g)_{p,a}\le \frac{c}{\nu^{2m}}\|\Delta_a^m g\|_{p,a},$$
where $c=c(n,a,m)$.
\end{lemma}

\proof
$$-\tilde{\Delta}_{y,a}^mg(x):=\sum_{l=1}^m(-1)^{l-1} \binom{m}{l} T^{\sqrt{l}y}_ag(x)-g(x)=\sum_{l=1}^md_lT^{\sqrt{l}y}_ag(x)-g(x),$$
where $\sum_{l=1}^m d_l=1$.\\
Let $\tilde{\eta}\in M(1, 1,a)$ be a nonnegative radial function such that $\|\tilde{\eta}\|_{1,a}=1$, moreover $\int_0^\infty \tilde{\eta}(\varrho)\varrho^{2m+|a|+n-1}d\varrho<\infty$. (E.g. $\tilde{\eta}(x)=C\left(\frac{\sin|x|^2}{|x|^2}\right)^M$, where $M\in \mathbb{N}$, $|a|+n+2m<2M$.) Let
$$Q_{\nu,m}(g)(x):=\int_{\mathbb{R}^n_+}\tilde{\eta}(y)\left(-\tilde{\Delta}_{\frac{y}{\nu},a}^mg(x)+g(x)\right)y^ady=\sum_{l=1}^m d_l\int_{\mathbb{R}^n_+}T^{\sqrt{l}\frac{y}{\nu}}_ag(x)\tilde{\eta}(y)y^ady$$ $$=\sum_{l=1}^m d_l\int_{\mathbb{R}^n_+}T^x_ag\left(\sqrt{l}\frac{y}{\nu}\right)\tilde{\eta}(y)y^ady=\int_{\mathbb{R}^n_+}\left(\sum_{l=1}^m d_l\left(\frac{\nu}{\sqrt{l}}\right)^{n+|a|}\tilde{\eta}\left(\frac{\nu}{\sqrt{l}}t\right)\right)T^x_ag(t)t^adt$$ $$=:K_{\nu,m,a}*_ag(x).$$
Now we estimate the difference of $g$ and $Q_\nu(g)$. By Minkowski's inequality we have
$$\|g-Q_{\nu,m}(g)\|_{p,a}=\left(\int_{\mathbb{R}^n_+}\left|\int_0^\infty \tilde{\eta}(\varrho)\int_{S^n_{1+}}\tilde{\Delta}_{\frac{\varrho \Theta}{\nu},a}^mg(x)\Theta^adS(\Theta)\varrho^{|a|+n-1}d\varrho\right|^px^adx\right)^{\frac{1}{p}}$$ $$\le |S^n_{1+}|_a\int_0^\infty \tilde{\eta}(\varrho)\varrho^{|a|+n-1}\|\tilde{\Delta}_{sph,\frac{\varrho \Theta}{\nu},a}^mg(x)\|_{p,a}d\varrho.$$
Thus by the assumption on $\tilde{\eta}$ and by Lemma \ref{l3} we have
$$\|g-Q_{\nu,m}(g)\|_{p,a}\le \frac{c}{\nu^{2m}}\int_0^\infty \tilde{\eta}(\varrho)\varrho^{|a|+n-1+2m}d\varrho \|\Delta_a^m g\|_{p,a}\le \frac{c}{\nu^{2m}}\|\Delta_a^m g\|_{p,a}.$$
Finally we have to show that $Q_{\nu,m}(g)\in M(\nu, p,a)$. $K_{\nu,m,a}\in E_e$ and so $Q_{\nu,m}(g)\in E_e$.
$$\|Q_{\nu,m}(g)\|_{p,a}\le \|K_{\nu,m,a}\|_{1,a}\|g\|_{p,a}.$$
Since
$$\|K_{\nu,m,a}\|_{1,a}\le \sum_{l=1}^m |d_l|\left(|S^n_{1+}|_a\int_0^\infty \left(\frac{\nu}{\sqrt{l}}\right)^{n+|a|}\tilde{\eta}\left(\frac{\nu}{\sqrt{l}}r\right)r^{n+|a|-1}dr\right)$$ $$=c(n,a,m)\|\tilde{\eta}\|_{1,a},$$
where
$$ c(n,a,m)<2^m|S^n_{1+}|_a.$$
Since $\tilde{\eta}\in M(1, 1,a)$, $\mathrm{supp}\hat{K}_{\nu,m,a}\subset B(0,\nu)$. Thus $\mathrm{supp}\hat{Q}_{\nu,m}(g)\subset B(0,\nu)$, which finishes the proof.

\medskip

\proof (of Lemma \ref{lbal}) In view of \eqref{t1}
$$\|\Delta_{sph,h,a}^mg\|_{p,a}\leq \|\Delta_{sph,h,a}^m(g-P_\nu(g)\|_{p,a}+\|\Delta_{sph,h,a}^mP_\nu(g)\|_{p,a}$$ $$\le 2^m\|g-P_\nu(g)\|_{p,a}+\|\Delta_{sph,h,a}^mP_\nu(g)\|_{p,a}.$$
According to \eqref{E}, Lemma \ref{l4} and Lemma \ref{l2}  if $h\in (0,t]$, $\nu=\frac{2}{h}$, we have
$$\|\Delta_{sph,h,a}^mg\|_{p,a}\leq 2^{2m}\frac{c}{\nu^{2m}}\|\Delta_a^m g\|_{p,a}+h^{2m}\|\Delta_a^m g\|_{p,a}\le c t^{2m}\|\Delta_a^m g\|_{p,a},$$
where $c\neq c(m).$ The proof is finished by taking supremum.

\medskip

To prove the right-hand side inequality in Theorem \ref{ok} we need the following lemmas.

\medskip

\begin{lemma}\label{l5} Let $f\in L^p_a$ $1\le p\le\infty$, $\nu>0$, $m\in\mathbb{N}$. Then we have
$$\|\Delta_a^mP_\nu(f)\|_{p,a}\le c(a,m,n) \nu^{2m}\omega_{sph,m}\left(f, \frac{1}{\nu}\right)_{p,a}.$$
\end{lemma}

\medskip

\begin{lemma}\label{l6} With the notation of the previous lemma we have
$$\|f-P_{\frac{\nu}{2}}(f)\|_{p,a}\le c(a,m,n)\|\Delta_{sph,\frac{1}{\nu},a}^mf\|_{p,a}.$$
\end{lemma}

\medskip

\proof (of the right-hand side of Theorem \ref{ok}) Since $P_\nu(f)\in M(\nu,p,a)$ (cf. \eqref{pnorm} and \eqref{Psupp}), by definition we have
$$ K\left(f, t^{2m}, L_a^p, W^{m,p}_{\Delta_a}\right)\le \|f-P_\nu(f)\|_{p,a}+t^{2m}\|\Delta_a^mP_\nu(f)\|_{p,a}.$$
In view of \eqref{E} and \eqref{Pnu} we have
$$\|f-P_\nu(f)\|_{p,a}\le c E_\nu(f)_{p,a}\le c\|f-P_{\frac{\nu}{2}}(f)\|_{p,a}.$$
Let $\nu=\frac{1}{t}$. Taking into consideration that $\|\Delta_{sph,\frac{1}{\nu},a}^mf\|_{p,a}\le \omega_{sph,m}\left(f, \frac{1}{\nu}\right)_{p,a}$, Lemma \ref{l5} and Lemma \ref{l6} imply
$$\|f-P_\nu(f)\|_{p,a}+t^{2m}\|\Delta_a^mP_\nu(f)\|_{p,a}$$ $$\le c\omega_{sph,m}\left(f, \frac{1}{\nu}\right)_{p,a}+ct^{2m}\nu^{2m}\omega_{sph,m}\left(f, \frac{1}{\nu}\right)_{p,a}=c\omega_{sph,m}\left(f, \frac{1}{\nu}\right)_{p,a}.$$

\medskip

\proof (of Lemma \ref{l5}) It is enough to prove that for all $h\in\left(0,\frac{1}{2\nu}\right]$,
\begin{equation}\label{DD}\|\Delta_a^mP_\nu(f)\|_{p,a}\le c \frac{1}{h^{2m}}\|\Delta_{sph,h,a}^mf\|_{p,a}.\end{equation}
Indeed, if \eqref{DD} is fulfilled, taking  $h=\frac{1}{2\nu}$, then taking supremum and considering the monotonicity of the modulus of smoothness, the lemma is verified.\\ To prove \eqref{DD} let take the Hankel transform of the left-hand side.  Since $h\le \frac{1}{2\nu}$, $\eta\left(\frac{|\xi|}{\nu}\right)=\eta\left(\frac{|\xi|}{\nu}\right)\eta(h|\xi|)$, thus we have
$$\mathcal{H}_a(\Delta_a^mP_\nu(f))=\frac{1}{h^{2m}}\eta\left(\frac{|\xi|}{\nu}\right)\frac{(-1)^m(h|\xi|)^{2m}\eta(h|\xi|)}{(1-j_{\frac{n+|a|}{2}-1}(h|\xi|))^m}(1-j_{\frac{n+|a|}{2}-1}(h|\xi|))^m\hat{f}(\xi).$$
Recall that
$$\mathcal{H}_a(\Delta_{sph,h,a}^mf)=(1-j_{\frac{n+|a|}{2}-1}(h|\xi|))^m\hat{f}(\xi).$$
Let us take into consideration that $$\frac{(-1)^m(h|\xi|)^{2m}\eta(h|\xi|)}{(1-j_{\frac{n+|a|}{2}-1}(h|\xi|))^m}\in \mathcal{S}_e,$$ and $$\left\|\mathcal{H}_a^{-1}\left(\frac{(-1)^m(h|\xi|)^{2m}\eta(h|\xi|)}{(1-j_{\frac{n+|a|}{2}-1}(h|\xi|))^m}\right)\right\|_{1,a}=\left\|\mathcal{H}_a^{-1}\left(\frac{(-1)^m(|\xi|)^{2m}\eta(|\xi|)}{(1-j_{\frac{n+|a|}{2}-1}(|\xi|))^m}\right)\right\|_{1,a}.$$
Since
$$\Delta_a^mP_\nu(f)=\frac{1}{h^{2m}}\varrho_\nu *_a\mathcal{H}_a^{-1}\left(\frac{(-1)^m(h|\xi|)^{2m}\eta(h|\xi|)}{(1-j_{\frac{n+|a|}{2}-1}(h|\xi|))^m}\right) *_a\Delta_{sph,h,a}^mf,$$
 by \eqref{Y} we obtain \eqref{DD}.

 \medskip

\proof (of Lemma \ref{l6})  We start again with Hankel transform.
$$\mathcal{H}_a\left(f-P_{\frac{\nu}{2}}(f)\right)=\frac{1-\eta\left(\frac{2|\xi|}{\nu}\right)}{\left(1-j_{\frac{n+|a|}{2}-1}\left(\frac{|\xi|}{\nu}\right)\right)^m}\left(1-j_{\frac{n+|a|}{2}-1}\left(\frac{|\xi|}{\nu}\right)\right)^m\hat{f}(\xi)$$ $$=\varphi\left(\frac{|\xi|}{\nu}\right)\mathcal{H}_a\left(\Delta_{sph,h,a}^mf\right).$$
$\varphi$ is decomposed by $\varphi_1+\varphi_2$, cf. \cite[(4.67)]{p}, where
$$\varphi_1=(1-\eta(2|\xi|))\left((1-j_{\frac{n+|a|}{2}-1}(|\xi|))^{-m}-Q_N(j_{\frac{n+|a|}{2}-1}(|\xi|))\right),$$ and
$$\varphi_2=(1-\eta(2|\xi|))Q_N(j_{\frac{n+|a|}{2}-1}(|\xi|)).$$
Here $Q_N(u)$ is a polynomial of degree $N$, and is defined as follows. With $N=N(n,a)$ large enough, $\frac{1}{(1-u)^m}=Q_N(u)+\sum_{j=N+1}^\infty b_ju^j$. In \cite[(4.68)-(4.70)]{p} is proved that
$$\left\|\mathcal{H}_{\frac{n+|a|}{2}-1}^{-1}\left(\varphi_1\left(\frac{r}{\nu}\right)\right)\right\|_{1,\frac{n+|a|}{2}-1,\mathbb{R}_+}=\left\|\mathcal{H}_{\frac{n+|a|}{2}-1}^{-1}\left(\varphi_1\left(r\right)\right)\right\|_{1,\frac{n+|a|}{2}-1,\mathbb{R}_+}$$
\begin{equation}\label{fi1n}=\left\|\mathcal{H}_{a}^{-1}\left(\varphi_1\left(\frac{|\xi|}{\nu}\right)\right)\right\|_{1,a,\mathbb{R}_+^n}=\left\|\mathcal{H}_{a}^{-1}\left(\varphi_1\left(|\xi|\right)\right)\right\|_{1,a,\mathbb{R}_+^n}<\infty,\end{equation}
i.e. the norm is finite and independent of $\nu$.\\
To treat $\varphi_2$ we introduce the spherical translation as follows.
\begin{equation}\label{str}T^r_{sph,a}f(x):=\frac{1}{|S^n_{1+}|_a}\int_{S^n_{1+}}T^{r\Theta}f(x)\Theta^a dS(\Theta).\end{equation}
It is clear that
$$\mathcal{H}_{a}((T^r_{sph,a})^kf)(\xi)=j_{\frac{n+|a|}{2}-1}^k(r|\xi|)\hat{f}(\xi).$$
Moreover in view of \eqref{t1} by Minkowski's inequality we have
\begin{equation}\label{tsn}\|T^r_{sph,a}f(x)\|_{p,a}\le \frac{1}{|S^n_{1+}|_a}\int_{S^n_{1+}}\|T^{r\Theta}f(x)\|_{p,a}\Theta^a dS(\Theta)\le \|f\|_{p,a}.\end{equation}
Thus for any $g\in\mathcal{S}'_e$
$$\varphi_2\left(\frac{|\xi|}{\nu}\right)\hat{g}(\xi)=\sum_{k=0}^N b_k j_{\frac{n+|a|}{2}-1}^k\left(\frac{|\xi|}{\nu}\right)\hat{g}(\xi)+\sum_{k=0}^N b_k \eta\left(\frac{2|\xi|}{\nu}\right)j_{\frac{n+|a|}{2}-1}^k\left(\frac{|\xi|}{\nu}\right)\hat{g}(\xi)$$ $$=\mathcal{H}_{a}\left(\sum_{k=0}^N b_k(T^r_{sph,a})^kg+\varrho_{\frac{\nu}{2}}*_a \sum_{k=0}^N b_k(T^r_{sph,a})^kg\right).$$
Finally we have
$$\|f-P_{\frac{\nu}{2}}(f)\|_{p,a}\le \left\|\mathcal{H}_{a}^{-1}\left(\varphi_1\left(\frac{|\xi|}{\nu}\right)\right)*_a\Delta_{sph,h,a}^mf\right\|_{p,a}$$ $$+\left\|\sum_{k=0}^N b_k(T^r_{sph,a})^k\Delta_{sph,h,a}^mf+\varrho_{\frac{\nu}{2}}*_a \sum_{k=0}^N b_k(T^r_{sph,a})^k\Delta_{sph,h,a}^mf\right\|_{p,a}.$$
Thus by \eqref{pnorm}, \eqref{fi1n}, \eqref{tsn} and \eqref{Y} we have
$$\|f-P_{\frac{\nu}{2}}(f)\|_{p,a}$$ $$\le \left(\|\mathcal{H}_{a}^{-1}(\varphi_1(|\xi|))\|_{1,a}+(1+\|\varrho_1\|_{1,a})\sum_{k=0}^N b_k\right)\|\Delta_{sph,h,a}^mf\|_{p,a} =c\|\Delta_{sph,h,a}^mf\|_{p,a}.$$

\medskip

Lemma \ref{l6} implies the next estimation.

\medskip

\begin{cor} Let $f\in L^p_a$ $1\le p\le\infty$, $\nu>0$, $m\in\mathbb{N}$. Then we have
$$ E_\nu(f)_{p,a}\le c \omega_{sph,m}\left(f, \frac{1}{\nu}\right)_{p,a}.$$
\end{cor}

\section{Bessel-Sobolev- and fractional Bessel-Sobolev spaces}

First, recalling Definition \ref{Sob} we make the following observations.

\noindent {\bf Remark.}

\noindent {\bf(1)} Take the next weighted Sobolev space:\\ $ W^{m,p}_{a}(\mathbb{R}^n_+):=\left\{f: \mathbb{R}^n_{+} \to \mathbb{R} : D^\alpha f \in  L^p_a, \ws |\alpha|\le m \right\}$. Notice, that the Bessel-Sobolev spaces are not coincide with the weighted Sobolev spaces. For instance, consider the fundamental solution of the Laplace-Besel operator, i.e. the solution to $\Delta_aE=\delta_a$, where $\delta_a$ stands for the weighted Dirac delta distribution, for which $\langle \delta_a,\varphi\rangle_a=\varphi(0)$. Then according to  \cite[Theorem 93]{ss}
$$E(x)=c(n,a)\left\{\begin{array}{ll}\ln |x|, \ws n+|a|=2\\ |x|^{2-n-|a|}, \ws n+|a|>2. \end{array}\right.$$
Let
$$f(x):=\left\{\begin{array}{ll}E(x), \ws |x|<1\\ 0, \ws |x|>2,\\ \end{array}\right.$$
and let $f$ be radial, $f\in C^2(\mathbb{R}^n_+)$.
Then $f, \Delta_a f \in L^p_a$ if $1<p<1+\frac{2}{n+|a|-2}$, $\frac{\partial}{\partial x_i}f \in L^p_a$ if $1\le p<1+\frac{1}{n+|a|-1}$ and $\frac{\partial^2}{\partial x_i\partial x_j}f \in L^p_a$ if $1\le p<1+\frac{1}{n+|a|}$. That is $f\in W^{1,p}_{\Delta_a}(\mathbb{R}^n_+)$, but $f\not\in W^{2,p}_{a}(\mathbb{R}^n_+)$ if $1+\frac{2}{n+|a|-2}\le p<1+\frac{1}{n+|a|}$. \\
On the other hand let $g(x)=e^{-\sum_{i=1}^n x_i}$. Then $g\in W^{2,p}_{a}(\mathbb{R}^n_+)$ for all $1\le p \le \infty$, but $g\not\in  W^{1,p}_{\Delta_a}(\mathbb{R}^n_+)$ if $p\ge n+|a|$.

\noindent {\bf(2)} Similarly to the classical case $W^{m,p}_{\Delta_a}(\mathbb({R}^n_+)$ are Banach spaces (c.f. \cite[Proposition 3.2]{sk}) which are separable if $1\le p<\infty$ and are reflexive and uniformly convex if $1<p<\infty$, cf. \cite[3.4]{a}.

\medskip

By interpolation of spaces (K-method, see \cite[(3), page 39]{bl}), $L^p_a$ and  $W^{1,p}_{\Delta_a}$ we arrive to the so-called Besov-type seminorm. If $p=q$, it is $$\left(\int_0^\infty\left(\frac{ K\left(f, u, L_a^p, W^{1,p}_{\Delta_a}\right)}{u^{\Theta}}\right)^p\frac{du}{u}\right)^{\frac{1}{p}}.$$
Replacing $u=t^2$ we get the following definition.\\
Let $0<s<2$.
$$\int_0^\infty\left(\frac{ K\left(f, t^{2}, L_a^p, W^{1,p}_{\Delta_a}\right)}{t^{s}}\right)^p\frac{dt}{t} \sim \int_0^\infty\left(\frac{ \omega_{s,1}(f,t)_{a,p}}{t^{s}}\right)^p\frac{dt}{t}$$ \begin{equation}=\int_0^\infty\left(\frac{\sup_{0<h<t}\|\Delta_{sph,h,a}f\|_{p,a} }{t^{s}}\right)^p\frac{dt}{t}=:\|f\|_{B_{\Delta_a,1,p,p}^s}^p.\end{equation}

We also define the spherical Gagliardo seminorm:
\begin{equation}\label{gags}[f]_{sph,s,p,a}:=\left(\int_0^\infty\left(\frac{\|\Delta_{sph,t,a}f\|_{p,a} }{t^{s}}\right)^p\frac{dt}{t}\right)^{\frac{1}{p}}\end{equation} $$=\frac{1}{|S^n_{1+}|_a}\left(\int_0^\infty\frac{\int_{\mathbb{R}^n_+}\left|\int_{S^n_{1+}}(f(x)-T^{t\Theta}_af(x))\Theta^adS(\Theta)\right|^px^adx}{t^{s p+1}}dt\right)^{\frac{1}{p}}$$ $$=\frac{1}{|S^n_{1+}|_a}\left(\int_{\mathbb{R}^n_+}\int_0^\infty\frac{\left|\int_{S^n_{1+}}(f(x)-T^{t\Theta}_af(x))\Theta^adS(\Theta)\right|^p}{t^{s p+n+|a|}}t^{n+|a|-1}dtx^adx\right)^{\frac{1}{p}}$$ \begin{equation}\label{gagw}\le\frac{1}{|S^n_{1+}|_a}\left(\int_{\mathbb{R}^n_+}\int_{\mathbb{R}^n_+}\frac{\left|f(x)-T^{y}_af(x)\right|^p}{|y|^{s p+n+|a|}}y^adyx^adx\right)^{\frac{1}{p}}=:[f]_{s,p,a},\end{equation}
which is the weighted Gagliardo seminorm. As a corollary of Theorem \ref{ok}, with the spherical Gagliardo seminorm we have the next equivalence.

\medskip

\begin{cor} If $f\in L^p_a$, $s\in (0,2)$, $$\|f\|_{B_{\Delta_a,1,p,p}^s}\sim [f]_{sph,s,p,a}.$$
\end{cor}

\proof
Of course, by Definition \ref{ms}, $[f]_{s,\gamma,p,a} \le  \|f\|_{B_{\Delta_a,1,p,p}^s}$. On the other hand in view of Definition \ref{k} we have
$$ \|f\|_{B_{\Delta_a,1,p,p}^s}^p \sim \int_0^\infty\frac{\left(\inf\{\|f-g\|_{p,a}+t^2\|\Delta_ag\|_{p,a} : g\in W^{1,p}_{\Delta_a}\}\right)^p}{t^{s p+1}}dt.$$
In view of \eqref{Pnu}
$$\|f\|_{B_{\Delta_a,1,p,p}^s}\le \int_0^\infty\frac{\left(\|f-P_{\frac{1}{2t}}(f)\|_{p,a}+t^2\|\Delta_a P_{\frac{1}{2t}}(f)\|_{p,a} \right)^p}{t^{s p+1}}dt.$$
According to Lemma \ref{l6} and \eqref{DD} we have
$$\|f\|_{B_{\Delta_a,1,p,p}^s}^p\le c\int_0^\infty\frac{\left(\|\Delta_{sph,t,a}f\|_{p,a}+\frac{t^2}{t^2}\|\Delta_{sph,t,a}f\|_{p,a}\right)^p}{t^{s p+1}}dt$$ $$=c(n,p,a) [f]_{sph,s,p,a}^p.$$

\medskip

Thus we define the fractional Bessel-Sobolev spaces as follows.

\begin{defi}Let $\frac{s}{2}\in(0,1)$, $1\le p <\infty$.
\begin{equation} W^{\frac{s}{2},p}_{\Delta_a}:=\left\{f\in L^p_a(\mathbb{R}^n_+) : \|f\|_{W^{\frac{s}{2},p}_{\Delta_a}}<\infty\right\},\end{equation}
where
\begin{equation} \|f\|_{W^{\frac{s}{2},p}_{\Delta_a}}=\left(\|f\|_{p,a}^p + [f]_{sph,s,p,a}^p\right)^{\frac{1}{p}}.\end{equation}
\end{defi}

\medskip

\noindent {\bf Remark.}

The so-called $K$-method ensures that  $W^{\frac{s}{2},p}_{\Delta_a}$ is also Banach space, cf. e.g. \cite[Theorem 3.4.2]{bl}.

\medskip

Below we need the fractional partial differential operators $(-\Delta_a)^{\frac{s}{2}}$ and $(I-\Delta_a)^{\frac{s}{2}}$, where $0<s<2$. Let $\varphi \in \mathcal{S}_e$.

\begin{equation}\label{fracl} (-\Delta_a)^{\frac{s}{2}}\varphi(x)=C(n,a,s)\int_{\mathbb{R}^n_{(+)}} \frac{\varphi(x)-T^y_a\varphi(x)}{|y|^{n+|a|+s}}y^ady, \end{equation}
see \cite{l} and in one dimension \cite{bg}.

In view of \cite[(46)]{dls} we have
\begin{equation}\label{Ds}(I-\Delta_a)^\frac{s}{2}\varphi(x)=\int_{\mathbb{R}^n_+}\frac{T^y_a\varphi(x)-\varphi(x)}{|y|^{n+|a|+s}}\omega_{a,-s}(|y|)y^ady +\varphi(x), \end{equation}
where
$$\omega_{a,\nu}(|y|)=c(n,a,\nu)K_{\frac{n+|a|-\nu}{2}}(|y|)|y|^{\frac{n+|a|-\nu}{2}},$$
$K_\nu$ is defined in \eqref{k0}.

The fractional operators above can also be defined using Hankel transformation. Let $\varphi \in \mathcal{S}_e$ again.
\begin{equation}\label{fraclh}\left(\mathcal{H}_a(-\Delta_a)^{\frac{s}{2}}\varphi\right)(\xi)=|\xi|^s\hat{\varphi}(\xi),\end{equation}
see \cite{l}, and
\begin{equation}\label{Dsh}\left(\mathcal{H}_a(I-\Delta_a)^\frac{s}{2}\varphi\right)(\xi)=(1+|\xi|^2)^{\frac{s}{2}}\hat{\varphi}(\xi),\end{equation}
see e.g. \cite{dls}.

Subsequently we need the next observation.

\medskip

\noindent {\bf Remark.}

Let $\varphi\in \mathcal{S}_e$, $h\in \mathcal{S}'_e$, $0\le s \le 2$ and let us denote by
$$D_s=(-\Delta_a)^{\frac{s}{2}} \ws \ws \mbox{or} \ws \ws D_s=(I-\Delta_a)^{\frac{s}{2}}.$$ Then we have
\begin{equation}\label{csere}\langle D_s\varphi, h\rangle_a=\langle \varphi, D_s h\rangle_a \end{equation}
and
\begin{equation}\label{csere2}D_s(\varphi*_a h)=D_s\varphi*_a h=\varphi*_a D_sh.\end{equation}

Indeed, let $\varphi, h\in \mathcal{S}_e$. Then \eqref{csere} follows from \eqref{fraclh} or \eqref{Dsh} and from the Parseval's formula, see \cite[page 20]{k}, in a standard way. After this we can naturally extend the operators to distributions by \eqref{csere}.\\
To derive \eqref{csere2} let $\varphi, h\in \mathcal{S}_e$ again. Taking into account \eqref{CH}, \eqref{fraclh} or \eqref{Dsh}, \eqref{csere2} can be proved by Hankel transform as well. Also taking into account \eqref{HT}, Hankel transform ensures that if $\varphi\in \mathcal{S}_e$, then $D_sT^t\varphi=T^tD_s\varphi$. Thus, in view of \eqref{csere} we have
$$\langle D_sh,T^t\varphi\rangle_a=\langle h, D_sT^t\varphi\rangle_a=\langle h,T^tD_s\varphi\rangle_a.$$
After the natural extension the first equality can be seen as follows. Let $\psi\in \mathcal{S}_e$.
$$\langle D_s(\varphi*_a h),\psi\rangle_a=\langle \varphi*_a h, D_s\psi\rangle_a=\langle h,\varphi*_aD_s\psi\rangle_a=\langle h,D_s\varphi*_a\psi\rangle_a=\langle D_s\varphi*_a h, \psi\rangle_a.$$

\medskip

\noindent The following theorem will be useful in what follows.

\begin{theorem}\label{tss} Let $0<\frac{s}{2}<1$; $1\le p<\infty$. Then $\mathcal{S}_e$ is a dense subset of $W^{1,p}_{\Delta_a}(\mathbb{R}^n_+)$ and $W^{\frac{s}{2},p}_{\Delta_a}(\mathbb{R}^n_+)$.
\end{theorem}

\proof
Obviously $\mathcal{S}_e\subset W^{1,p}_{\Delta_a}$, and independently, Lemma \ref{lbal} ensures that $\mathcal{S}_e\subset W^{\frac{s}{2},p}_{\Delta_a}$.

First we prove that $C_0^\infty(\mathbb{R}^n_+)$ is dense in $W^{1,p}_{\Delta_a}(\mathbb{R}^n_+)$. Following the standard chain of ideas we use Bessel convolution instead of the standard one. Let $\tilde{\varepsilon}>0$ be arbitrary, $u\in W^{1,p}_{\Delta_a}(\mathbb{R}^n_+)$. Then there are $\delta$ and $R$ (depending on $u$ and $\tilde{\varepsilon}$) such that $\|u\|_{W^{1,p}_{\Delta_a},\mathbb{R}^n_+\setminus \Omega_{\tilde{\varepsilon}}}<\tilde{\varepsilon}$, where $\Omega_{\tilde{\varepsilon}}=\mathbb{R}^n_+\cap B(0,R)\setminus\{x \in \mathbb{R}^n_+ : \mathrm{dist}(x,\partial \mathbb{R}^n_+)\le \delta\}$. Let us cover $\overline{\Omega_{\tilde{\varepsilon}}}$ with balls of radius at most $\frac{\delta}{2}$. Select a finite covering $\Omega_{\tilde{\varepsilon}}\subset \cup_{j=1}^mB_j$ and let $\{\Psi_j\}_{j=1}^m$ be the corresponding partition of unity. Define $u_j:=\Psi_j u$, and $\eta_{a,\varepsilon}(x):=\frac{k_a}{\varepsilon^{n+|a|}}\eta\left(\frac{x}{\varepsilon}\right)$, where $\eta(x)=k_ae^{-\frac{1}{1-|x|^2}}$, if $|x|<1$, $\eta(x)=0$ if $|x|\ge 1$ and $\int_{\mathbb{R}^n_+}\eta(x)x^a=1$. Let $\varepsilon<\frac{\delta}{2}$ and for sake of abbreviation let $d(x,t,\vartheta)=\sum_{i=1}^nx_i^2+t_i^2-2x_it_i\cos\vartheta_i$. Considering that if $\varepsilon<\|x-t\|\le \sqrt{d(x,t,\vartheta)}$, then $T^t\eta_{a,\varepsilon}(x)=0$, $\mathrm{supp}\eta_{a,\varepsilon}(x)*_au_j\subset \mathbb{R}^n_+$, so $\eta_{a,\varepsilon}(x)*_au_j\in C_0^\infty(\mathbb{R}^n_+)$.\\
Let $\varphi_j\in C_0^\infty(\mathbb{R}^n_+)$ such that $\|u_j-\varphi_j\|_{p,a}\le \epsilon_j$. Then $$\|u_j-\eta_{a,\varepsilon}(x)*_au_j\|_{p,a}\le \|u_j-\varphi_j\|_{p,a}\left(1+\|\eta_{a,\varepsilon}\|_{1,a}\right)+\|\varphi_j-\eta_{a,\varepsilon}*_a\varphi_j\|_{p,a}.$$
Since $\|\eta_{a,\varepsilon}\|_{1,a}=\|\eta_{a,1}\|_{1,a}$, i.e. is independent of $\varepsilon$, it is enough to deal with the second term above. Taking into account the norm of the Bessel translation we have
$$\left|\varphi_j-\eta_{a,\varepsilon}*_a\varphi_j(x)\right|=\left|\int_{\mathbb{R}^n_+}(\varphi_j(x)-\varphi_j(t))T^t\eta_{a,\varepsilon}(x)t^adt\right|\le \sup_{|x-t|<\varepsilon}|\varphi(x)-\varphi(t)|,$$
thus
$$\|\varphi_j-\eta_{a,\varepsilon}*_a\varphi_j\|_{p,a}\le 2\epsilon_j |\mathrm{supp}\varphi_j|_a.$$
In view of \eqref{csere2} we have $\Delta_a(\eta_{a,\varepsilon}(x)*_au_j)=\Delta_a\eta_{a,\varepsilon}(x)*_au_j=\eta_{a,\varepsilon}(x)*_a\Delta_a u_j$, in distribution sense.
Since $u\in W^{1,p}_{\Delta_a}$, $\Delta_au_j\in L^p_a$. Indeed, recalling that for all $\varphi\in \mathcal{S}_e$, $u\in\mathcal{S}'_e$ $\langle \Delta_au, \varphi\rangle_a=\langle u, \Delta_a\varphi\rangle_a$,  we have\\ $\|\Delta_a u\|_{p,a}=\sup_{\varphi\in \mathcal{S}_e, \|\varphi\|_{p',a}\le 1}\int_{\mathbb{R}^n_+}|u(x)\Delta_a\varphi(x)|x^adx=:K$. Thus we get
$$\left|\int_{\mathbb{R}^n_+}\Delta_a u_j \varphi x^adx\right|=\left|\int_{\mathbb{R}^n_+} u_j \Delta_a\varphi x^adx\right|\le \|\psi_j\|_\infty K.$$
That is we can apply the previous part to $\Delta_a u_j$. Since the partition is finite, choosing $\tilde{\varepsilon}$, $\varepsilon$, $\epsilon$ small enough, $\sum_{j=1}^m \eta_{a,\varepsilon}(x)*_au_j$ is can be arbitrarily close to $u$ in $W^{1,p}_{\Delta_a}(\mathbb{R}^n_+)$-norm.

Let $0<s<2$. According to \cite[Lemma 7.18]{a} it is enough to define an appropriate sequence of linear operators. Let $D_k=\mathbb{R}^n_+\cap B(0,k)\setminus\{x\in \mathbb{R}^n_+: \mathrm{dist}(x,\partial \mathbb{R}^n_+)\le\frac{1}{k}\}$. With the above procedure we construct $\varphi_k:=\sum_{j=1}^{m_k}\eta_{a,\varepsilon_k}*_au_j\in C_0^\infty(\mathbb{R}^n_+)$, such that $\|u-\varphi_k\|_{W^{1,p}_{\Delta_a}, D_k}<\frac{1}{k}$. Since $\|u\|_{W^{1,p}_{\Delta_a},\mathbb{R}^n_+\setminus D_k}\to 0$ when $k\to \infty$, with $P_ku:=\varphi_k$ we have that $\|P_ku-u\|_{B_i}\to 0$, where $B_1=L^p_a(\mathbb{R}^n_+)$, $B_2=W^{1,p}_{\Delta_a}(\mathbb{R}^n_+)$, so $\mathcal{S}_e\subset B_1\cap B_2$ is dense in $W^{\frac{s}{2},p}_{\Delta_a}(\mathbb{R}^n_+)$, $0<\frac{s}{2}<1$.

\section{Bessel B-potential spaces}

To introduce Bessel-potential spaces we define the appropriate kernel function $G_{a,\nu}$. The Bessel kernel is
\begin{equation}\label{g1}G_{a,\nu}(x):=\frac{2^{\frac{n-|a|-\nu}{2}+1}}{\Gamma\left(\frac{\nu}{2}\right)\prod_{i=1}^n\Gamma(\alpha_i+1)}\frac{K_{\frac{n+|a|-\nu}{2}}(|x|)}{|x|^{\frac{n+|a|-\nu}{2}}},\end{equation}
where $K_\nu$ is defined by \eqref{k0}.

\eqref{k1} and \eqref{k2} imply that if $\nu>0$, $G_{a,\nu}\in L^1_a(\mathbb{R}^n_+)$. Moreover
\begin{equation}\label{GH}G_{a,\nu}(x)=\mathcal{H}_a^{-1}\left(\frac{1}{(1+|\xi|^2)^{\frac{\nu}{2}}}\right),\end{equation}
see e.g. \cite[(5)]{dls}.

\medskip

Now we are in position to define the next spaces.

\medskip

\begin{defi}
\begin{equation}\label{panu} L^{\nu,p}_{a}:=L^{\nu,p}_{a}(\mathbb{R}^n_+)=\{f :f=G_{a,\nu}*_a g; \ws   g\in L^p_a\}, \ws \ws \ws \|f\|_{p,a,\nu}:=\|g\|_{p,a}.\end{equation}
\end{defi}

\noindent {\bf Remark}

Direct calculations show that if $g \in C_0^\infty$, then  $G_{a,s}*_ag \in \mathcal{S}_e$.

\medskip

\noindent For integers the connection between Bessel-Sobolev and Bessel B-potential spaces is described by the following statement, see \cite[Lemma 2.]{h} or \cite[Theorem 4.4]{sk}.

\begin{statement}\label{l1}\cite[Lemma 2.]{h} Let $a_i>0$ $i=1, \dots , n$ and $m$ be a positive integer, $1<p<\infty$. Then (with equivalent norms)
$$W^{m,p}_{\Delta_a}=L^{2m,p}_{a}.$$
\end{statement}

\medskip

The non-integer case is a little different, cf. \cite[Theorem 7.63]{a}.

\medskip

\begin{theorem}\label{t22} Let $1<p<\infty$, $s>\varepsilon>0$ and $0<s<2$ be arbitrary. Then
$$L_a^{s+\varepsilon,p} \hookrightarrow W^{\frac{s}{2},p}_{\Delta_a} \hookrightarrow L_a^{s-\varepsilon,p}.$$
\end{theorem}

\medskip

Together with Theorem \ref{tss} the following proposition implies Theorem \ref{t22}.

\begin{proposition}\label{p2} Let $1<p<\infty$, $s>\varepsilon>0$, $0<s<2$ be arbitrary.\\
Let $f\in  L^{s+\varepsilon,p}_{a}$. Then
\begin{equation}\label{WL}[f]_{sph,s,p,a}\le c \|f\|_{p,a,s+\varepsilon}.\end{equation}
Lef $f\in\mathcal{S}_e$. Then there is a $g\in L^p_a$ such that $f=G_{a,s-\varepsilon}*_ag$, and we have
\begin{equation}\label{LW}\|g\|_{p,a}\le c \|f\|_{W^{\frac{s}{2},p}_{\Delta_a}}.\end{equation}
In the inequalities above the constants are independent of $f$ and $g$.
\end{proposition}

\medskip

To prove the above proposition we need the following lemma, see \cite{dls} or \cite{n}.

\medskip

\begin{lemma}\cite[Theorem 3]{dls}\label{ldls} For all $\varphi\in \mathcal{S}_e$, $s>0$ we have
\begin{equation}\label{inv}(1-\Delta_a)^{\frac{s}{2}}(G_{a,s}*_a\varphi)=G_{a,s}*_a(1-\Delta_a)^{\frac{s}{2}}\varphi=\varphi.\end{equation}
\end{lemma}

\medskip

\proof (of Proposition \ref{p2}) First we prove \eqref{WL}. Let $f=G_{a, s+\varepsilon}*_a g$. In view of Theorem \ref{tss} and the above remark we can assume that $g\in C^\infty_0(\mathbb{R}^n_+)$. Then by \eqref{CH}, \eqref{spht} and \eqref{GH} we have
$$\mathcal{H}_a\left(\Delta_{sph,r,a}f\right)(\xi)=\left(1-j_{\frac{n+|a|}{2}-1}(r|\xi|\right)\frac{1}{(1+|\xi|^2)^{\frac{s+\varepsilon}{2}}}\hat{g}(\xi)$$ $$=\mathcal{H}_a\left((\Delta_{sph,r,a}G_{a,s+\varepsilon})*_ag\right)(\xi).$$
$$\|\Delta_{sph,r,a}f\|_{p,a}\le \|\Delta_{sph,r,a}G_{a,s+\varepsilon}\|_{1,a}\|g\|_{p,a},$$
and according to \eqref{gags} we have
$$[f]_{sph,s,p,a}\le \|g\|_{p,a}\left(\int_0^\infty\frac{1}{r^{sp+1}}\|\Delta_{sph,r,a}G_{a,s+\varepsilon}\|_{1,a}^p dr\right)^{\frac{1}{p}}=:\|g\|_{p,a}H(G,a,s,\varepsilon,p).$$
$$H(G,a,s,\varepsilon,p)$$ $$\le \left(\int_0^1\frac{1}{r^{sp+1}}\|\Delta_{sph,r,a}G_{a,s+\varepsilon}\|_{1,a}^p dr\right)^{\frac{1}{p}}+\left(\int_1^\infty\frac{1}{r^{sp+1}}\|\Delta_{sph,r,a}G_{a,s+\varepsilon}\|_{1,a}^p dr\right)^{\frac{1}{p}}$$ $$=H_1+H_2.$$
By \eqref{t1}
$$\|\Delta_{sph,r,a}G_{a,s+\varepsilon}\|_{1,a}\le \|G_{a,s+\varepsilon}\|_{1,a}<\infty,$$
so
$$H_2\le c.$$
Denoting by $t:=|x|$ let us recall that by \eqref{besatl} and \eqref{rh} we have
$$\mathcal{H}_a(\Delta_{sph,r,a}G_{a,s+\varepsilon})$$ $$=c(n,a)\mathcal{H}_{\frac{n+|a|}{2}-1}\left(\frac{K_{\frac{n+|a|-(s+\varepsilon)}{2}}(t)}{t^{\frac{n+|a|-(s+\varepsilon)}{2}}}-T^r_{\frac{n+|a|}{2}-1}\frac{K_{\frac{n+|a|-(s+\varepsilon)}{2}}(t)}{t^{\frac{n+|a|-(s+\varepsilon)}{2}}}\right).$$
Since $G_{a,s+\varepsilon}\in L^1_a(\mathbb{R}^n_+)$ and $G_b(t):=\frac{K_{\frac{n+|a|-(s+\varepsilon)}{2}}(t)}{t^{\frac{n+|a|-(s+\varepsilon)}{2}}}\in L^1_{n+|a|-1}(\mathbb{R}_+)$, by \eqref{bi} and \cite[Lemma 3.2]{esk} we have
$$\|\Delta_{sph,r,a}G_{a,s+\varepsilon}(|x|)\|_{1,a}=c_{a,n}$$ $$\times\int_0^\infty\left|\int_0^\pi \left(G_{b}(t)-G_{b}\left(\sqrt{t^2+r^2-2rt\cos\vartheta}\right)\right)\sin^{n+|a|-1}\vartheta d\vartheta\right| t^{n+|a|-1}dt$$ $$=c_{a,n}\left(\int_0^{\frac{r}{2}}|\cdot| t^{n+|a|-1}dt+\int_{\frac{r}{2}}^{2r} |\cdot| t^{n+|a|-1}dt+\int_{2r}^\infty|\cdot| t^{n+|a|-1}dt\right)$$ $$=h_1(r)+h_2(r)+h_3(r).$$
$$H_1\le  \left(\int_0^1\left(\frac{h_1(r)}{r^{s+\frac{1}{p}}}\right)^p dr\right)^{\frac{1}{p}}+ \left(\int_0^1\left(\frac{h_2(r)}{r^{s+\frac{1}{p}}}\right)^p dr\right)^{\frac{1}{p}}+ \left(\int_0^1\left(\frac{h_3(r)}{r^{s+\frac{1}{p}}}\right)^p dr\right)^{\frac{1}{p}}.$$
Let $u:=\sqrt{t^2+r^2-2rt\cos\vartheta}$. Since $|r-t|\le u \le r+t$, in the first integral $u>\frac{r}{2}$. In view of \eqref{k1}
$$\int_0^1\left(\frac{h_1(r)}{r^{s+\frac{1}{p}}}\right)^p dr$$ $$\le c(n,a,s)\int_0^1\frac{1}{r^{sp+1}}\left(\int_0^{\frac{r}{2}}\left(\frac{1}{t^{n+|a|-s-\varepsilon}}-\frac{1}{r^{n+|a|-s-\varepsilon}}\right)t^{n+|a|-1}dt\right)^pdr$$ $$\le  c(n,a,s)\int_0^1\frac{r^{(s+\varepsilon)p}}{r^{sp+1}}dr=c(n,a,s,\varepsilon)<\infty.$$
According to \eqref{k3} we have
$$\int_0^1\left(\frac{h_2(r)}{r^{s+\frac{1}{p}}}\right)^p dr$$ $$\le c\int_0^1\frac{1}{r^{sp+1}}\left(r \max\left\{\frac{1}{r^{n+|a|-s-\varepsilon}},\right.\right.$$ $$\left.\left. \int_0^\pi G_{b}\left(\sqrt{t^2+r^2-2rt\cos\vartheta}\right)\sin^{n+|a|-1}\vartheta d\vartheta\right\} r^{n+|a|-1}\right)^pdr.$$
Let $t=(1+\delta)r$, $-\frac{1}{2}\le \delta \le 1$. Recalling the notation above, $u^2=r^2\left(\delta^2+4(1+\delta)\sin^2\frac{\vartheta}{2}\right)\ge r^2\sin^2\vartheta$. Thus by \eqref{k1} we get
$$\int_0^\pi G_{b}\left(\sqrt{t^2+r^2-2rt\cos\vartheta}\right)\sin^{n+|a|-1}\vartheta d\vartheta \le c \frac{1}{r^{n+|a|-s-\varepsilon}}\int_0^\pi\frac{\vartheta ^{n+|a|-1}}{\vartheta^{n+|a|-s-\varepsilon}}d\vartheta$$ $$=c(s,\varepsilon)\frac{1}{r^{n+|a|-s-\varepsilon}},$$
and so as above we have
$$\int_0^1\left(\frac{h_2(r)}{r^{s+\frac{1}{p}}}\right)^p dr\le  c\int_0^1\frac{1}{r^{sp+1}}\left(r \frac{1}{r^{n+|a|-s-\varepsilon}}r^{n+|a|-1}\right)^pdr=\int_0^1r^{\varepsilon p-1}dr<\infty.$$
Finally, in view of \eqref{k2} we have
$$\int_0^1\left(\frac{h_3(r)}{r^{s+\frac{1}{p}}}\right)^p dr \le c \int_0^1\frac{1}{r^{sp+1}}\left(\int_{2r}^\infty \frac{e^{-t}}{t^{1-s-\varepsilon}}dt\right)^pdr\le c\int_0^1r^{\varepsilon p-1}dr<\infty.$$
Thus
$$H_1<c,$$
and \eqref{WL} is proven.

To prove \eqref{LW} we use \eqref{Ds}.\\
Recalling that $f\in\mathcal{S}_e$, according to Lemma \ref{ldls} we have
\begin{equation}\label{Ginv}f=G_{a, s-\varepsilon}*_a(I-\Delta_a)^{\frac{s-\varepsilon}{2}}f.\end{equation}
Thus we have to estimate $\|(I-\Delta_a)^{\frac{s-\varepsilon}{2}}f\|_{p,a}$. In view of \eqref{Ds} we need the following estimation.
$$\left(\int_{\mathbb{R}^n_+}\left|\int_{\mathbb{R}^n_+}\frac{T^y_af(x)-f(x)}{|y|^{n+|a|+s-\varepsilon}}\omega_{a,-(s-\varepsilon)}(|y|)y^ady\right|^px^adx\right)^{\frac{1}{p}}$$ $$\le \left\|\int_{|y|\le 1\cap \mathbb{R}^n_+}(\cdot)\right\|_{p,a}+\left\|\int_{|y|> 1\cap \mathbb{R}^n_+}(\cdot)\right\|_{p,a}=I+II.$$
\eqref{k2} and \eqref{t1} ensure that
$$II\le c \|f\|_{p,a}.$$
In view of \eqref{k1} we have
$$I\le\left\|\int_0^1\int_{S_{1+}^n}(T^{r\theta}f(x)-f(x))\theta^adS(\theta)\frac{1}{r^{1+s-\varepsilon}}dr\right\|_{p,a}$$ $$\le c(n,a)\left\|\int_0^1\frac{|\Delta_{s,r,a}f(x)|}{r^{1+s-\varepsilon}}dr\right\|_{p,a}.$$
Since $(\cdot)^p$ is convex, we have
$$\left(\int_0^1\frac{|\Delta_{s,r,a}f(x)|}{r^{s+\frac{\varepsilon}{p-1}}r^{1-\frac{p\varepsilon}{p-1}}}dr\right)^p\le \left(\int_0^1\frac{1}{r^{1-\frac{p\varepsilon}{p-1}}}\right)^{p-1}\int_0^1\frac{|\Delta_{s,r,a}f(x)|^p}{r^{p\left(s+\frac{\varepsilon}{p-1}\right)}r^{1-\frac{p\varepsilon}{p-1}}}dr,$$
cf. \cite[(6.14.1)]{hlp}. That is
$$I\le c(n,a,p,\varepsilon)\left(\int_{\mathbb{R}^n_+}\int_0^1\frac{|\Delta_{s,r,a}f(x)|^p}{r^{ps+1}}drx^adx\right)^{\frac{1}{p}}=c\left(\int_0^1\frac{\|\Delta_{s,r,a}f(x)\|_{p,a}^p}{r^{ps+1}}dr\right)^{\frac{1}{p}},$$
so the estimations of $I$ and $II$ together imply \eqref{LW}.

\medskip

\proof(of Theorem \ref{t22}.) Since $s>0$, by \eqref{k1} and \eqref{k2} $\|f\|_{p,a}\le c \|f\|_{p,a,S+\varepsilon}$, i.e. if $f\in\mathcal{S}_e$,
\begin{equation}\label{nn} \|f\|_{p,a,s-\varepsilon}\le c \|f\|_{W^{\frac{s}{2},p}_{\Delta_a}} \le c\|f\|_{p,a,s+\varepsilon}.\end{equation}
It remains only to prove that $W^{\frac{s}{2},p}_{\Delta_a} \subset L_a^{s-\varepsilon,p}$.\\
According to theorem \ref{tss} and the remark after Definition \ref{panu} $\mathcal{S}_e$ is dense in both Banach spaces (with respect to the appropriate norm).
If Let $f\in W^{\frac{s}{2},p}_{\Delta_a}$. $\{f_n\}\subset \mathcal{S}_e$ such that $f_n \to f$ in $W^{\frac{s}{2},p}_{\Delta_a}$-norm, then it is also Cauchy in
$L_a^{s-\varepsilon,p}$-norm, i.e. there is an $\tilde{f}:=\lim_{n\to\infty} f_n$ in $L_a^{s-\varepsilon,p}$. Since $f,\tilde{f}\in L^p_a$ we have $\|f-\tilde{f}\|_{p,a}\le \|f-f_n\|_{p,a}+\|f_n-f_m\|_{p,a}+\|\tilde{f}-f_m\|_{p,a}\le c(\|f-f_n\|_{W^{\frac{s}{2},p}_{\Delta_a}}+\|f_n-f_m\|_{p,a,s-\varepsilon}+\|\tilde{f}-f_m\|_{p,a,s-\varepsilon}$. Thus $f=\tilde{f}$ a.e..

\section{Removable sets}

\begin{defi} We say that a compact set $K\subset \mathbb{R}^n_+$ is removable with respect to the partial differential operator $L$ in $L^p_a(\mathbb{R}^n_+)$ if for any $u\in L^p_a(\mathbb{R}^n_+\setminus K)$, $Lu=0$ on $\mathbb{R}^n_+\setminus K$ there is a $\tilde{u}\in L^p_a(\mathbb{R}^n_+)$ such that $\tilde{u}|_{\mathbb{R}^n_+\setminus K}=u$ and $Lu=0$ on $\mathbb{R}^n_+$.
\end{defi}

\medskip

To characterize removable sets the next definition of different capacities is useful, cf. \cite[Definitions 2.2.6 and 2.7.1]{ah}. Recall that $\mathcal{S}_e\subset L^{p,s}_a$ is dense.

\begin{defi}Let $K\subset \mathbb{R}^n_+$ be compact. $U\subset \mathbb{R}^n_+$ stands for a neighbourhood of $K$.
$$C_{p,a,s}(K):=\inf\left\{\|f\|_{p,a,s}^p: f \in \mathcal{S}_e ; \ws f(x)\ge 1, \ws \forall x \in K\right\}.$$

$$N_{p,a,s}(K):=\inf\left\{\|f\|_{p,a,s}^p: f\in \mathcal{S}_e; \ws \exists \ws U,\ws K\subset U, \ws f(x)= 1, \ws \forall x\in U\right\}$$
\end{defi}

\medskip

\noindent {\bf Remark.}

\noindent - It can be readily seen that if $N_{p,a,s}(K)=0$, then $\lambda(K)=0$, see e.g. \cite{h}.\\
- According to \cite[Theorem 2.5.1]{ah} we have
\begin{equation}\label{cmu}C_{p,a,s}(K)^{\frac{1}{p}}=\sup\left\{\mu(K): \|G_{a,s}*_a \mu\|_{p',a}\le 1\right\},\end{equation}
where $\mu$ is a positive measure supported on $K$. \\
- Obviously, $C_{p,a,s}(K)\le N_{p,a,s}(K)$.\\
- If $C_{p,a,s}(K)=0$, then $N_{p,a,s}(K)=0$ as well, see \cite{h}. Indeed, let $h\in C^\infty(\mathbb{R}_+)$,\\$ 0\le h\le 1$, such that
$$h(t)=\left\{\begin{array}{ll} 0, \ws 0<t\le \frac{1}{2},\\ 1, \ws t\ge 1\end{array}\right.$$ If $C_{p,a,s}(K)=0$, then for all $\varepsilon>0$ there is a nonnegative $f\in C_0^\infty$ such that $G_{a,s}*_af\ge 1$ on $K$, and $\|f\|_{p,a}^p\le \frac{\varepsilon}{2}$. Then we can assume that there is a neighbourhood $U\supset K$, such that $G_{a,s}*_af\ge 1$ on $U$, and $\|f\|_{p,a}^p\le \varepsilon$. Then $h(G_{a,s}*_af)$ satisfies the requirements of the definition of $N_{p,a,s}(K)$. Moreover, since
$$\|G_{a,s}*_af\|_{p,a}^p\ge\int_{\{G_{a,s}*_af\ge \frac{1}{2}\}}\frac{1}{2^p}x^adx,$$
in view of \eqref{k2} we have
$$c(s,a,p)2^p\varepsilon\ge 2^p\|G_{a,s}*_af\|_{p,a}^p\ge \int_{\{G_{a,s}*_af\ge \frac{1}{2}\}}x^adx\ge\int_{\mathbb{R}^n_+}(h(G_{a,s}*_af(x))^px^adx,$$
which implies that $N_{p,a,s}(K)=0$.

\medskip

\note

Let $0<s_k<2$ and $a_k$ be real numbers ($k=0, \dots , m$). We define the next fractional partial differential operator as follows.
\begin{equation}\label{Lu} Lu:=\sum_{k=1}^m a_k(I-\Delta_a)^{\frac{s_k}{2}}u.\end{equation}

\medskip

\begin{theorem} Let $L$ be as in \eqref{Lu} and $K\subset \mathbb{R}^n_+$ be compact. Let $2>s>\max_{k=1}^m s_k$. If $C_{p',a,s}(K)=0$, then $K$ is removable   with respect to $L$ $L^p_a(\mathbb{R}^n_+)$. If $m=1$, then $K$ is removable (with respect to $L=(I-\Delta_a)^{\frac{s}{2}}$) if and only if $C_{p',a,s}(K)=0$.
\end{theorem}

\proof
Let $C_{p',a,s}(K)=0$. Then the last remark ensures that $N_{p',a,s}(K)=0$, i.e. for all $\varepsilon>0$ there is a $\varphi=G_{a,s}*_a g\in \mathcal{S}_e$ such that $\|\varphi\|_{p,a,s}<\varepsilon$ and $\varphi\equiv 1$ on some neighbourhood $U$ of $K$, cf. \eqref{panu}. Let $O$ be an open set such that $K\subset U\subset O\subset \mathbb{R}^n_+$, and let $\psi\in C_0^\infty(O)$. Then $(1-\varphi)\psi\in  C_0^\infty(O\setminus K)$. Let $u\in\mathcal{S}'_e\cap L^p_a(O\setminus K)$ such that $Lu=0$ on $O\setminus K$. In view of \eqref{csere} we have
$$\langle Lu, (1-\varphi)\psi\rangle_a=\langle u, L(1-\varphi)\psi\rangle_a=\langle u, L\psi\rangle_a-\langle u, L\varphi\psi\rangle_a=0.$$
Thus by the remark above we have
$$\left|\langle u, L\psi\rangle_a\right|=\left|\langle u, L\varphi\psi\rangle_a\right|\le \|u\|_{p,a,O}\|L\varphi\psi\|_{p',a}.$$
By linearity it is enough to deal with a single term. Since
$$(I-\Delta_a)^{\frac{s_k}{2}}\varphi\psi=(I-\Delta_a)^{\frac{s_k}{2}}((G_{a, s}*_ag)\psi),$$
by \eqref{csere} we compute
$$\sup_{\tau\in \mathcal{S}_e, \|\tau\|_{p,a}\le 1}\langle G_{a, s}*_ag, \psi(I-\Delta_a)^{\frac{s_k}{2}}\tau \rangle_a$$ $$ \le \|\psi\|_\infty \sup_{\tau\in \mathcal{S}_e, \|\tau\|_{p,a}\le 1}\langle |G_{a, s}*_ag|, |(I-\Delta_a)^{\frac{s_k}{2}}\tau |\rangle_a$$ $$ =c\sup_{\tau\in \mathcal{S}_e, \|\tau\|_{p,a}\le 1}\langle (I-\Delta_a)^{\frac{s_k}{2}}(G_{a, s}*_ag), \tau \rangle_a $$ $$=c\sup_{\tau\in \mathcal{S}_e, \|\tau\|_{p,a}\le 1}\langle G_{a, s-s_k}*_ag, \tau \rangle_a\le c \|g\|_{p',a}\le c\varepsilon.$$
The last but one step can be seen by Hankel transform again. Indeed, if $g\in C_0^\infty$, the inversion formula, \eqref{bi} can be applied on both sides, moreover by \eqref{GH}, \eqref{Dsh} and \eqref{CH} $\mathcal{H}_a\left((I-\Delta_a)^{\frac{s_k}{2}}(G_{a, s}*_ag)\right)(\xi)=\mathcal{H}_a\left( G_{a, s-s_k}*_ag\right)(\xi)=\frac{1}{(1+|\xi|^2)^{\frac{s-s_k}{2}}}\hat{g}(\xi)$. Since $s-s_k>0$, in view of \eqref{k1}, $ G_{a, s-s_k}\in L^1_a(\mathbb{R}^n_+)$.

Let $m=1$. Assume that $C_{p',a,s}(K)>0$. In view of \eqref{cmu} there is a positive measure $\mu$ supported on $K$ such that $G_{a, s}*_a\mu\in L^p_a(\mathbb{R}^n_+)$. Together with \eqref{cca} and \eqref{csere2}, Lemma \ref{ldls} implies that the fundamental solution of the operator $(1-\Delta_a)^{\frac{s}{2}}$ is $G_{a,s}\in \mathcal{S}'_e$. Then $L(G_{a, s}*_a\mu)$ is supported on $K$, i.e. $K$ is not removable.

\medskip

\vspace{5mm}

{\small{\noindent Both authors:\\
Department of Analysis and Operations Research,\\ Institute of Mathematics,\newline
Budapest University of Technology and Economics \newline
 M\H uegyetem rkp. 3., H-1111 Budapest, Hungary.

 \vspace{3mm}

\noindent \'A. P. Horv\'ath: g.horvath.agota@renyi.hu\\
M. Chegaar: mouna.chegaar@edu.bme.hu}

\end{document}